\documentclass[12pt]{article}

\usepackage{geometry}
\usepackage{color}
\usepackage{amssymb}

\newif\ifforward
\forwardtrue
\forwardfalse

\title{{\normalsize\tt\hfill\jobname.tex}\\
On SPDE and backward filtering equations for SDE systems 
(direct approach)
\\ ~ \\
A. Yu. Veretennikov\footnote{
University of Leeds, UK; National Research University Higher School of Economics and Institute for Information Transmission Problems, Moscow, Russia, 
e-mail: a.veretennikov @ leeds.ac.uk
}
}

\begin{document}
\maketitle

\newcommand{\new}{\newcommand}
\new{\beq}{\begin{equation}}
\new{\eeq}{\end{equation}}
\new{\beqr}{\begin{eqnarray}}
\new{\eeqr}{\end{eqnarray}}
\new{\Fw}{{\cal F}^{\tilde w}}
\new{\rhoi}{\rho^{-1}_{t_{i},T}}

\newtheorem{Th}{Theorem}
\newtheorem{Le}{Lemma}

\newtheorem{Theorem}{Theorem}
\newtheorem{Lemma}{Lemma}

\begin{abstract}
A direct approach to linear backward filtering equations 
for SDE systems is proposed. This preprint is a corrected version of the paper 1995 in the LMS Lecture Notes combined with another paper by the author on the direct approach to linear SPDEs for SDEs. 
\end{abstract}

\section{Introduction}

Filtering theory is one of the main sourses of stochastic partial differential equations (SPDE's).
In this paper the filtering problem for the simplest 
two-dimensional stochastic differential equation system is considered, 

\beqr\label{eq1}
dX_t=f(X_t)dt+dw^1_t,\quad X_0=x,
\nonumber\\ \\ \nonumber
dY_t=h(X_t)dt+dw^2_t,\quad Y_0=y,
\eeqr
where functions $f$ and $h$ are smooth and bounded, $w^1$ and $w^2$ are 
independent standard Wiener processes on some probability space 
$(\Omega, {\cal F}, ({\cal F}_t,\, t\geq 0), \mathbb{P})$; initial data $X_0=x$ and $Y_0=y$ are assumed non-random. (In fact, both may be distributed being mutually independent with $(w^1, w^2)$ and for the component $Y$ this may be helpful in filtering, but we do not pursue this goal here.) The problem is to describe 
the estimate of unobservable process $X_t$ via observable component 
$Y_t,\, 0\leq s\leq t$, which is optimal in the mean--square sense, 
i.e., 
$$
m_t\equiv m_{0,t}=\mathbb{E}[g(X_t)|{\cal F}^Y_t],\quad t\geq 0.
$$
In fact, the answer is known even 
for more general situations, see \cite{Ro}. In particular, $m_t$ may be
represented via backward stochastic differential equations, which makes sense if we are interested in an optimal estimation for 
some fixed time $t$; in this case we should find the solution 
of our backward SPDE and substitute there the trajectory of our observation process.

In this paper 
we present a {\em direct} approach to such a representation, using a similar idea 
for an equation without filtering, that is, for a completely observed SDE trajectory. This preprint is an improved version of the paper \cite{Ve95} presented along with the main lines of the calculus from \cite{Ve}. The matter is that the standard way -- as in \cite{KR, P, Ro} -- is to write down the SPDE, then establish existence and uniqueness of solution in  appropriate (Sobolev) classes, then apply Ito's (or Ito--Wentzell's) formula and, hence, justify that this solution, indeed, coincides with the desired conditional expectation. Apparently, this way assumes that somehow the equation should be known in advance. What the direct approach provides is exactly how to derive the equation ``by hand'' without reference to any big theory.
Note that there is a paper \cite{KZ} with a very similar title; yet, this is a different {\em direct approach,} which also stems from Krylov's idea of representing solutions of SDEs as solutions of linear SPDEs, see \cite{Kry}, \cite{Ro}, \cite{Ve}.

The paper consists of four sections. Number one is the Introduction; the second one contains the main result about filtering SPDEs as well as two auxiliary Lemmata; the third one is devoted to the proof of the Lemma \ref{Th2} (the second Lemma is a well known result with a reference provided), and the fourth one contains the proof of the main result -- the Theorem \ref{Th1}.

\section{Main result and auxiliary lemmata}

Due to Girsanov's theorem, process $Y_t,\,0\leq t\leq T$ is a Wiener process
on some the probability space with some new measure $(\Omega, {\cal F}, ({\cal F}_t, 0\leq t\leq T),\tilde{\mathbb{P}})$ 
(see below).

\begin{Th}[backward SPDE]\label{Th1}
Let $f,h \in C^3_b$. Then 
the process $m_t$ may be represented as follows:
\beq\label{eq2}
m_T=\frac{v^g(0,x)}{v^1(0,x)},
\eeq
where the processes  $v^g$ and $v^1$ satisfy the following linear
backward stochastic differential equation (the same for both functions):
\beq\label{eq3}
-dv^g(t,x)=\left[\frac12 v^g_{xx}(t,x)+f(x)v^g_x(t,x)\right] dt+
h(x)v^g_{}(t,x)\star dY_t,\quad 0\leq t\leq T,
\eeq
with initial data
\beq\label{eq4}
v^g(T,x)=g(x),\quad x\in \mathbb{R}^1. 
\eeq
\end{Th}
Note that the denominator in (\ref{eq2}) is strictly positive a.s. as a conditional expectation of a strictly positive random variable with respect to some new probability measure. This will be commented in the proof.

\noindent
Here in (\ref{eq3}) $\int \cdot \star dY_t$ means ``backward'' stochastic Ito integral,
i.e., a normal ``regular'' stochastic Ito integral with inverse time, see 
\cite{Kry, Ro}. It may be formally defined, for example, by the formula 
\begin{eqnarray}\label{inverseint}
\int_0^T h(x)v^g_{}(t,x)\star dY_t := 
\int_0^T h(x)\tilde v^g_{}(t,x) d\tilde Y_{t}, 
 \nonumber \\\\ \nonumber
\tilde Y_{t}=Y_T-Y_{T-t}, \quad \tilde v^g_{}(t,x)=v^g_{}(T-t,x), 
\end{eqnarray}
where $\int_0^T h(x)\tilde v^g_{}(t,x) d\tilde Y_{t}$ is a standard It\^o's integral. (The only small nuance is that this integral {\em might} be naturally defined up to the $\pm$ sign -- which relates simply to how a Wiener process in the inverse time is defined --  and, clearly,  this sign would also affect the sign in the last term of the equation (\ref{eq3}); this will be commented later.)
The function $v^1$ has its terminal condition $v^1(T,x)\equiv 1$ and satisfies
the same SPDE (\ref{eq3}). Notice that the random function $v^g(t,x)$ is,
in fact, $F^{w^{1},w^{2}}_{t,T}$-adapted ({\em not} $F^{w^{1},w^{2}}_{0,t}$-adapted); therefore, the integral above makes
sense exactly as a classical standard It\^o's one (cf. \cite[Theorem 6.3.1]{Ro}).

~

Before the proof we recall another Krylov and Rozovsky's result -- the Lemma \ref{Th2} below -- concerning
multidimensional SDEs (see \cite{Kry}, \cite{Ro}, \cite{Ve}). 

Let $(Z^{s,z}_t,\, t\geq s,\, s\geq 0,\,
z\in R^d)$ be the family of $d$-dimensional processes depending on the parameters $(s,z)$ and  satisfying the following multidimensional SDEs:

\beq\label{eq5}
dZ^{s,z}_t=b(Z^{s,z}_t)dt+\sigma(Z^{s,z}_t)dw_t,\quad t\geq s,\quad
Z^{s,z}_s=z,
\eeq
where $b$ is a bounded smooth $d$-dimensional vector, $\sigma$ is a matrix
$d\times d_1$, $w_t$ is a $d_1$-dimensional Wiener process, $d, d_1\ge 1$; there are neither any other restrictions on the values $d$ and $d_1$, nor any non-degenerabilty condition is assumed.
We will use the following different notations for the same value:
\[
Z^{s,z}_t\equiv Z(s,t,z),
\]
and for $t=T$ also
\[
Z^{s,z}_T=u(s,z).
\]
Recall that here $T$ is fixed throughout the text, and that the multidimensional 
setting is essential: 
we will need it in the proof of the Theorem \ref{Th1} with $d=2, \, d_1=1$. 

\begin{Le}\label{Th2}
Let $b,\sigma \in C^3_b$.
Then the random field \(Z^{s,z}_T\) is continuous in all variables \((s,T,z)\). Moreover, 
continuous partial derivatives exist, the gradient vector $\partial_z Z^{s,z}_{t}=: Z_z(s,t,z)$ and the Hessian matrix $\partial^2_{zz} Z^{s,z}_{t} =: Z_{zz}(s,t,z)$, 
and the process $u(s,z)$ satisfies an SPDE
\beqr\label{eq6}
&  \displaystyle -du(t,z)=\left[\frac12 (\sigma\sigma^*)_{ij}(z)u_{z_{i}z_{j}}(t,z)
+b^i(z)u_{z_{i}}(t,z)\right]dt
\nonumber\\ \\ \nonumber
& +\sigma_{ij}(z)u_{z_{i}}(t,z)\star dw^j_t,
\eeqr
with a terminal condition
\beq\label{eq7}
u(T,z)\equiv z.
\eeq
Here $\sigma^*$ means the matrix $\sigma$ transposed, and the 
equation (\ref{eq6}) holds true for each component of the vector $u(t,z) = (u^1(t,z), \ldots, u^d(t,z))^*$.
\end{Le}

The {\em direct} approach to this result may be found in \cite{Ve}, and the main lines of its  proof will be recalled below for the convenience of the reader.

~

Further, we will use the Bayes
representation for conditional expectations, 
also known as Kallianpur--Striebel's formula, see \cite{Ro}.

\begin{Le}\label{Le1}
Let the Borel functions \(h, g\) be bounded. Then the following representation is valid a.s.:
\[
m_T=\frac{\tilde \mathbb E[g(X_T)\rho^{-1}|{\cal F}^Y_T]}{\tilde \mathbb E[\rho^{-1}|{\cal F}^Y_T]},
\]
where $\tilde \mathbb E$ is the expectation with respect to the measure $\tilde  \mathbb P$:
$d\tilde  \mathbb P =\rho  \mathbb dP$, with
\[
\rho\equiv\rho_{0T}=\exp\left(-\int_0^Th(X_t)dw^2_t - \frac12 \int_0^T|h(X_t)|^2dt\right).
\]
\end{Le}

Recall that due to Girsanov's theorem the process $(Y_t,\, 0\leq t\leq T)$
is a Wiener process on probability space $(\Omega, {\cal F}, ({\cal F}_t,\, 0\leq t\leq T),\tilde \mathbb P)$, independent of $w^1$ and, in general, with a non-zero starting value.
Denote $\tilde w_t=Y_t-y$. Then on the space $(\Omega, {\cal F},({\cal F}_t,\,
0\leq t\leq T), \tilde \mathbb P)$ our system (\ref{eq1}) has the form

\beqr\label{eq8}
& dX_t=f(X_t)dt+dw^1_t, \quad X_0=x,
\nonumber\\ \\ \nonumber
& dY_t=d\tilde w_t,\quad Y_0=y,
\eeqr
with two independent Wiener processes  \((w^1, \tilde w)\).

\section{Direct proof of Lemma \ref{Th2}}
The terminal condition (\ref{eq7}) is straightforward: \(u(T,z) = Z^{T,z}_T=z\).

~

\noindent
Further, due to \cite{BlFr} (or \cite{Kunita} under relaxed assumptions) the random field $Z^{s,z}_{t}$ is continuous in $(s,t,z)$ for $s\le t$ and $z\in R^d$ and, moreover, it admits classical continuous partial derivatives $\partial_z Z^{s,z}_{t}=: Z_z(s,t,z)$ and $\partial^2_{zz} Z^{s,z}_{t} =: Z_{zz}(s,t,z)$.

~

\noindent
Note that continuity in all arguments is required, in particular, for the justification of a substitution of $Z^{0,z}_s$ into the initial value of $Z^{s,\cdot}_t$, that is, 
\begin{equation}\label{evolution}
Z^{0,z}_t = Z^{s, Z^{0,z}_s}_t.
\end{equation}
However, the equation contains two derivatives in the state variable, so we assumed the conditions, which guarantee the existence of these two (also classical) derivatives. The equation 
(\ref{evolution}) itself follows easily from the uniqueness of solution of the related SDEs and, indeed, from the Markov property, which  is a standard feature of any SDE with a unique solution in strong or weak sense (see, e.g., \cite{Kry-MP}).

~

\noindent
Further, what is stated in the Theorem, by definition 
may be written in the integral form  with $a = \sigma \sigma^*/2$ as follows (recall about the minus sign in the left hand side of the equation (\ref{eq6})), 
\begin{eqnarray}\label{eq8}
Z(t,T,z) - Z(T,T,z) 
= \int_t^T Z_z(s,T,z)\sigma(x)\star dW_s 
  \nonumber \\ \\ \nonumber
+\int_t^T \left[Z_z(s,T,z)b(z) + \mbox{Tr} \,Z_{zz}(s,T,z) a(z)\right] ds.
\end{eqnarray}

~

\noindent
In the sequel the general case $d, d_1\ge 1$ is presented; we will use it in the case $d=2, \, d_1=1$ in the proof of the Theorem \ref{Th1}.  
To show (\ref{eq8}), let us split the interval $[0,t]$ by small partitions $t = t_0 < t_1 < \dots < t_{n+1}= T$, and let us write down the identity, 
$$
Z(t,T,z) - Z(T,T,z) 
= \sum_{i=0}^{n} \left(Z(t_i,T,z) - Z(t_{i+1},T,z)\right),
$$
and consider each term after substituting  
$$
Z(t_i,T,z) = Z(t_{i+1},T,Z(z,t_i,t_{i+1})),
$$
and using Hadamard's form of Newton--Leibnitz' formula (also known as the First Theorem of the Calculus),
\begin{eqnarray*}
& \displaystyle F(\Delta)=F(0) + \int_0^1 \nabla F(\alpha  \Delta) \Delta d\alpha \hspace{1cm}
 \\ \\
&  \displaystyle =F(0) +   \int_0^1 \left(\nabla F(0) \Delta + \alpha  \int_0^1 \Delta^* \left(\mbox{div} \nabla F_{}(\alpha\beta)\right) \Delta)   d\beta \right) d\alpha, 
\end{eqnarray*}
where \(\nabla F(\alpha  \Delta) \Delta\) is the inner product of the two vectors, and 
\(\Delta^* \left(\mbox{div} \nabla F_{}(\alpha\beta)\right) \Delta\) is the Hessian matrix of \(F\) multiplied by the {\em vector} \(\Delta\) on the right and by the transposed \(\Delta^*\) on the left. 
Hence, we write, 
\begin{eqnarray*}
&Z(t_{i+1},T,Z(z,t_i,t_{i+1})) - Z(t_{i+1},T,z) = Z_z(t_{i+1},T,z)z_i 
 \\ \\
&  \displaystyle + \int_0^1\int_0^1 \alpha z_i^*  Z_{zz}(t_{i+1},T,z+ \alpha\beta z_i) z_i d\alpha  d\beta, 
\end{eqnarray*}
or, in the coordinate notations,
\begin{eqnarray*}
&Z(t_{i+1},T,Z(z,t_i,t_{i+1})) - Z(t_{i+1},T,z) = Z_{z^k}(t_{i+1},T,z)z_i^k 
 \\ \\
&  \displaystyle + \int_0^1\int_0^1 \alpha z_i^k  Z_{z^k z^\ell}(t_{i+1},T,z+ \alpha\beta z_i) z_i^\ell d\alpha  d\beta, 
\end{eqnarray*}
where summation over repeated indices is assumed (Einstein's convention),
$$
z_i := Z(t_i,t_{i+1},z) - z,
$$
and the equation is understood component-wise, i.e., for each component of the vector \(Z\). 
Denote also
$$
\tilde z_i := \sigma(z)(W_{t_{i+1}} - W_{t_{i}}) + b(z) (t_{i+1} - t_{i}),
$$
and let $\Delta W_{t_i} = W_{t_{i+1}} - W_{t_{i}}$.
By virtue of standard estimates in stochastic analysis it follows, 
$$
\sup_i \mathbb E|\tilde z_i - z_i|^2 \le \frac{C}{n^2}. 
$$
Hence, we get 
\begin{eqnarray}\label{eqzi}
Z(t_{i+1},T,Z(z,t_i,t_{i+1})) - Z(t_{i+1},T,z) = Z_{z^k}(t_{i+1},T,z)\tilde z_i^k 
 \nonumber \\ \\ \nonumber
+ \int_0^1\int_0^1 \alpha Z_{z^k z^\ell}(t_{i+1},T,x+ \alpha\beta z_i)\tilde z_i^k \tilde z_i^\ell d\alpha  d\beta 
+ o(1/n),  
\end{eqnarray}
where $o(1/n)$ is understood in the square mean sense. We have, $Z_z(t_{i+1},T,z)\tilde z_i \approx Z_z(t_{i+1},T,z)\sigma(z) (\Delta W_{t_i} + b(z)\Delta t_i)$; and 
$Z_{zz}(t_{i+1},T,z)\tilde z_i^2  \approx Z_{zz}(t_{i+1},T,z)\sigma^2(z) (\Delta W_{t_i})^2$. In all cases the sign $``\approx \ldots''$ means $``= \ldots + o(\max_i \Delta t_i)''$ with $o(\max_i \Delta t_i)$ in the square mean sense as $\max_i \Delta t_i \to 0$.

Recall the definition of the backward integral for $\xi_t \in {\cal F}^W_{t,T}$: 
$$
\int_0^T \xi(t)\star dW(t) := \int_0^T \xi_T(s) dW_T(s), 
$$
where 
$
\xi_T(s) = \xi(T-s), \; W_T(s) = W(T) - W(T-s). 
$
So, the integral approximations for the right hand side integral here with $0=t_0 < t_1 < \cdots < t_n = T$ read, 
\begin{eqnarray*}
&\sum_i \xi_T (t_i) (W_T(t_{i+1}) - W_T(t_i)) 
 \\\\
&= \sum_i \xi_{} (T-t_i) (W(T)-W_{}(T-t_{i+1}) - W(T)+W_{}(T-t_i))
 \\\\
&= \sum_i \xi_{} (T-t_i) (W_{}(T-t_i)-W_{}(T-t_{i+1}))
 \\\\
&= \sum_i \xi(t'_i)(W(t'_{i}) - W(t'_{i+1})),
\end{eqnarray*}
where $t'_i = T-t_i$. Note that this may be used as a simplified definition of stochastic integral, at least, for continuous $\xi(t)$. Since $0=t'_n < t'_{n-1} < \cdots < t'_0 = T$, the right way to understand integral approximations in terms of original processes in direct time is $\displaystyle \sum_i \xi(t_{i+1})(W(t_{i+1}) - W(t_{i})) 
=  \sum_i \xi(t_{i+1}) \Delta W_{t_i}$ (recall that $\Delta W_{t_{i}} = W_{t_{i+1}} - W_{t_{i}}$).
So, after summation over $i$ in (\ref{eqzi}), we obtain $Z(t,T,z) - Z(T,T,z)$ in the left hand side, and the following three terms (all component-wise) in the right hand side,
$$
\sum_i Z_{z^k}(t_{i+1},T,z)\sigma^{k\ell}(z)\Delta W^\ell_{t_{i}} \stackrel{sq.mean}{\to}    \int_t^T Z_z(s,T,z)\sigma(z)\star d W_s, 
$$
$$
\sum_i Z_{z^k}(t_{i+1},T,z)b^k(z)\Delta t_{i} \stackrel{sq.mean}{\to} 
\int_t^T Z_z(s,T,z)b(z) ds, 
$$
and
$$
\frac12\sum_i Z_{z^k z^\ell}(t_{i+1},T,z)\sigma^{kj}(z)\sigma^{\ell j}(z)(\Delta W^j_{t_{i}})^2 \stackrel{sq.mean}{\to} 
\int_t^T \mbox{Tr} \left(Z_{zz}(s,T,z)a(z)\right) ds,
$$
as $\max_i \Delta t_i \to 0$. Here \(\displaystyle \frac12\) is due to \(\displaystyle \int_0^1 \alpha \, d\alpha = \frac12\). 
So, we obtain (\ref{eq8}), as required.
The Lemma \ref{Th2} is proved.

\section{Direct proof of Theorem \ref{Th1}}

\noindent
{\bf 1}. Denote

\[
v^g(s,x) = \tilde \mathbb  E[g(X^{s,x}_T)\rho^{-1}_{s,T}|{\cal F}^Y_{s,T}].
\]
Then, 
\[
v^g(T,x) = \tilde  \mathbb   E[g(X^{T,x}_T)\rho^{-1}_{T,T}|{\cal F}                             ^Y_{T,T}] = g(x).
\]
In fact, what we want to establish is exactly the following equality (for each $T>0$ and any $0 \le t_0\le T$):

\beqr\label{eq89}
& v^g(t_0,x)-v^g(T,x) 
 \nonumber\\ \\ \nonumber
& \displaystyle = \int_{t_0}^T \left[\frac12 v^g_{xx}(t,x)+f(x)v^g_x(t,x)\right]dt
+ \int_{t_0}^T h(x) v^g(t,x)\star d\tilde w_t.
\eeqr
Let us use the identity

\[
v^g(t_0,x)-v^g(T,x)=\sum_{i=1}^N (v^g(t_{i-1},x)-v^g(t_i,x)),
\]
for any partition $t_0<t_1<\ldots <t_N=T$.
Consider one term from this sum: we have,

\begin{eqnarray*}
v^g(t_{i-1},x)-v^g(t_i,x) \hspace{3cm}
 \\\\
= \tilde  \mathbb  E[\rho^{-1}_{t_{i-1},T} g(X(t_{i-1},T,x))|{\cal F}^Y_{t_{i-1},T}] -
 \tilde  \mathbb  E[\rho^{-1}_{t_{i},T} g(X(t_i,T,x))|{\cal F}^Y_{t_{i},T}]
 \\\\
=\tilde  \mathbb  E[\rho^{-1}_{t_{i-1},T} g(X(t_{i-1},T,x))|{\cal F}^{\tilde w}_{t_{i-1},T}]-
\tilde  \mathbb  E[\rho^{-1}_{t_{i},T} g(X(t_i,T,x))|{\cal F}^{\tilde w}_{t_{i},T}].
\end{eqnarray*}

~

\noindent
{\bf 2}. Using continuity of the family $X(s,T,x)$ with respect to all
variables and existence of two continuous partial derivatives  with respect to $x$ (see \cite{BlFr}) we get a.s. by virtue of Taylor's expansion,

\begin{eqnarray*}
&X(t_{i-1},T,x)=X(t_i,T,X^{t_{i-1},x}_{t_{i}})
\\\\
&= X^{t_{\color{red}i},x}_T+X_x(t_i,T,x)(X^{t_{i-1},x}_{t_{i}}-x)
\\\\
& \displaystyle + \frac12 X_{xx}(t_i,T,x)(X^{t_{i-1}, x}_{t_{i}}-x)^2 + \alpha^1_i
\\\\
&=X^{t_{i},x}_T + X_x(t_i,T,x)(f(x)\Delta t_i + \Delta w^1_{t_{i}})
\\\\
& \displaystyle + \frac12  X_{xx}(t_i,T,x)\Delta t_i + \alpha^2_i,
\end{eqnarray*}
where $\Delta t_i=t_i-t_{i-1}$, $\Delta w^j_{t_{i}}=w^j_{t_{i}}-
w^j_{t_{i-1}}$, and $|\alpha^1_i|+|\alpha^2_i|=o(\Delta t_i)$ in the
mean-square sense. 
Hence, 
\begin{eqnarray*}
&g(X(t_{i-1},T,x))=g(X(t_i,T,X^{t_{i-1},x}_{t_{i}}))
\\\\
&= g\left(X^{t_{i-1},x}_T+X_x(t_i,T,x)(X^{t_{i-1},x}_{t_{i}}-x) \right.
 \\\\
&  \displaystyle \left. +\frac12 X_{xx}(t_i,T,x)(X^{t_{i-1},x}_{t_{i}}-x)^2 + \alpha^1_i\right)
 \\\\
&=g\left(X^{t_{i},x}_T + X_x(t_i,T,x)(f(x)\Delta t_i + \Delta w^1_{t_{i}}) \right.
 \\\\
& \displaystyle \left. + \frac12  X_{xx}(t_i,T,x)\Delta t_i + \alpha^2_i\right)
  \\\\
&= g(X^{t_{i},x}_T)  
 \\ \\
& \displaystyle + g_x (X^{t_{i},x}_T) \left(X_x(t_i,T,x)(f(x)\Delta t_i + \Delta w^1_{t_{i}})
+ \frac12  X_{xx}(t_i,T,x)\Delta t_i\right)
 \\\\
& \displaystyle + \frac12 g_{xx}(X^{t_{i},x}_T)
(\Delta w^1_{t_{i}})^2 
+ o(\Delta t_i) 
 \\\\
&= g(X^{t_{i},x}_T)  
 \\ \\
& \displaystyle + g_x(X^{t_{i},x}_T) \left(X_x(t_i,T,x)(f(x)\Delta t_i + \Delta w^1_{t_{i}})
+ \frac12  X_{xx}(t_i,T,x)\Delta t_i\right)
 \\\\
& \displaystyle +  \frac12  g_{xx}(X^{t_{i},x}_T) \Delta t_i  
+ o(\Delta t_i). 
\end{eqnarray*}
Denote $V(s,t,x) = g(X(s,t,x))$. Then,   assuming that $g\in C^2$, we have,  
$$
V_x = g_x X_x; \quad V_{xx} = g_x X_{xx} + g_{xx}X_x^2, 
$$
where we dropped the arguments in $g_{x}$, $g_{xx}$, $X_x$, and $X_{xx}$ for brevity. 
So, 
\begin{eqnarray*}
&g(X(t_{i-1},T,x))=
V(t_{i-1},T,x) 
\\\\
&= g(X^{t_{i},x}_T)  
 \\ \\
& \displaystyle + g_x(X^{t_{i},x}_T) (X_x(t_i,T,x)(f(x)\Delta t_i + \Delta w^1_{t_{i}})
+ \frac12  X_{xx}(t_i,T,x)\Delta t_i
 \\\\
& \displaystyle +  \frac12(X_x(t_i,T,x))^2 g_{xx}(X^{t_{i},x}_T) \Delta t_i + o(\Delta t_i) 
 \\\\
& \displaystyle = V(t_i,T,x) + V_x(t_i,T,x)(f(x)\Delta t_i + \Delta w^1_{t_{i}}) + \frac12 V_{xx}(t_i,T,x) \Delta t_i + o(\Delta t_i).
\end{eqnarray*}
Here and earlier $o(\Delta t_i)$ is understood in the mean square sense. 
The obtained relation means that the conditional expectation for $V = g(X)$ should satisfy the same SPDE as for $X$ itself, just with another terminal condition. 

~

\noindent
{\bf 3}. Thus, 

\begin{eqnarray*}
& \tilde  \mathbb  E[\rho^{-1}_{t_{i-1},T} g(X (t_{i-1},T,x))|{\cal F}^{\tilde w}_{t_{i-1},T}]
\\\\
& =\tilde  \mathbb  E [\rho^{-1}_{t_{i-1},T} g(X(t_i,T,X^{t_{i-1},x}_{t_{i}}))|
{\cal F}^{\tilde w}_{t_{i-1},T}] 
= \tilde  \mathbb  E [\rho^{-1}_{t_{i-1},T} V(t_{i-1},T,x) |
{\cal F}^{\tilde w}_{t_{i-1},T}]  
\\\\
& = \tilde  \mathbb  E [\rho^{-1}_{t_{i-1},T}
\{V(t_i,T,x)+V_x(t_i,T,x)f(x)\Delta t_i
 \\\\
&  \displaystyle +V_x((t_i,T,x)\Delta w^1_{t_{i}} + \frac12 V_{xx}(t_i,T,x)\Delta t_i\}
|{\cal F}^{\tilde w}_{t_{i-1},T}]
+ \alpha^3_i.
\end{eqnarray*}
Here again, $\alpha^3_i=o(\Delta t_i)$ in the mean square sense, i.e., $(\mathbb E |\alpha^3_i|^2)^{1/2}=o(\Delta t_i)$.

~

\noindent
{\bf 4}. Now, we would like to replace $\rho^{-1}_{t_{i-1},T}$ by
$\rho^{-1}_{t_{i},T}$. For this aim we apply the Lemma \ref{Th2} to 
the process $(X^{s,x}_{t},\,\rho^{-1}_{s,t},\,t\geq s)$. More precisely, 
let us note that this two-dimensional process satisfies the following SDE 
system: 

\beqr\label{eq9} 
& dX^{s,x}_t=f(X^{s,x}_t)dt+dw^1_t,\quad X^{s,x}_s=x,
\nonumber\\ \\
& d\rho^{-1}_{s,t}=h(X^{s,x}_t)\rho^{-1}_{s,t}d\tilde w_t,  \quad \rho^{-1}_{s,s}=1, 
\nonumber
\eeqr
with $s\leq t\leq T$. 
Indeed, $\rho^{-1}_{s,t}$ has the following representation:

\[
\rho^{-1}_{s,t}=\exp\left(\int_s^t h(X^{s,x}_r)d\tilde w_r 
- \frac12  \int_s^t |h(X^{s,x}_r)|^2 dr\right). 
\]
Let us consider a bit more general set of processes
$\{(X^{s,x}_t,\rho^{-1,\xi}_{s,t})\}$ which satisfy SDE's

\beqr\label{eq10}
& dX^{s,x}_t=f(X^{s,x}_t)dt+dw_t,\quad X^{s,x}_s=x,
\nonumber\\ \\ \nonumber
& d\rho^{-1,\xi}_{s,t}=h(X^{s,x}_t)\rho^{-1,\xi}_{s,t}d\tilde w_t, \quad \rho^{-1,\xi}_{s,s}=\xi,
\eeqr
for $s\leq t\leq T$, with $\xi>0$. In fact, $X^{s,x}_t$ here is the same as earlier, and 
$\rho^{-1,\xi}_t$ has the following representation: 

\[
\rho^{-1,\xi}_{s,t}=
\xi\exp\left(\int_s^th(X^{s,x}_r)d\tilde w_r - 
\frac12 \int_s^t|h(X^{s,x}_r)|^2\,dr\right)=\xi\rho^{-1}_{s,t}.
\] 
Then due to the Lemma \ref{Th2}, we get 

\new{\rhoxi}{\rho^{-1,\xi}}
\begin{eqnarray*}
& \displaystyle -d_s\rho^{-1,\xi}_{s,t} = \left[\frac12 h^2(x)(\rho^{-1,\xi}_{s,t})^2
(\rho^{-1,\xi}_{s,t})_{\xi\xi} + \frac12 (\rho^{-1,\xi}_{s,t})_{xx}\right.
 \\\\
& \displaystyle \left. + f(x)(\rho^{-1,\xi}_{s,t})_x\right]dt + (\rhoxi_{s,t})_x\star dw^1_t
+h(x) \xi (\rhoxi_{s,t})_\xi\star d\tilde w_t.
\end{eqnarray*}
Note that, in fact, $(\rhoxi_{s,t})_\xi=\rho^{-1}_{s,t}$ and
$(\rhoxi_{s,t})_{\xi\xi}=0$.
Hence,

\begin{eqnarray*}
& \displaystyle  -d_s\rhoxi_{s,t}=\left[\frac12 (\rhoxi_{s,t})_{xx}+f(x)(\rhoxi_{s,t})_x\right]dt 
\\\\
& +(\rhoxi_{s,t})\star dw^1_t+h(x)\xi \rho^{-1}_{s,t}\star d\tilde w_t.
\end{eqnarray*}
So, we get,

\begin{eqnarray*}
&\rhoxi_{t_{i-1},T}-\rhoxi_{t_{i},T}\equiv -\Delta\rhoxi_{t_{i}}
\\\\
& \displaystyle = \left[\frac12 (\rhoxi_{t_{i},T})_{xx}+f(x)(\rhoxi_{t_{i},T})_x\right]\Delta t_i
\\\\
&+ \rhoxi_{t_{i},T}\Delta w^1_{t_{i}} + h(x) \xi \rhoxi_{t_{i},T}\Delta
\tilde w_{t_{i}} + \alpha^4_i,
\end{eqnarray*}
with a similar $o(\Delta t_i)$ property for $\alpha^4_i$ as for previous $\alpha^1_i, \alpha^2_i, \alpha^3_i$.
Below we will use this assertion with $\xi=1$, that is, 
\begin{eqnarray*}
&\rho^{-1}_{t_{i-1},T}-\rhoi\equiv -\Delta\rhoi
\\\\
& \displaystyle = \left[\frac12 (\rhoi)_{xx}+f(x)(\rhoi)_x\right]\Delta t_i
\\\\
&+ \rhoi\Delta w^1_{t_{i}} + h(x)\rhoi\Delta
\tilde w_{t_{i}} + \alpha^4_i,
\end{eqnarray*}

~

\noindent
{\bf 5}. Now, we obtain

\begin{eqnarray*}
\tilde  \mathbb  E [\rho^{-1}_{t_{i},T} V(t_{i-1},T,x)|{\cal F}^{\tilde w}_{t_{i-1},T}] \hspace{4cm} 
\\\\
=\tilde  \mathbb  E\left[\left\{V(t_{i},T,x) +
(f(x)V_x(t_i,T,x) + \frac12 V_{xx}(t_i,T,x))\Delta t_i
+V_x(t_i,T,x)\Delta w^1_{t_{i}}\right\}\times \right.
\\\\
\left.\times
\left\{\rhoi + \left(\frac12 (\rhoi)_{xx} + f(x)(\rhoi)_x\right) \Delta t_i \hspace{3cm}\right.\right.
\\\\
\left.\left. \phantom{\frac12}+(\rhoi)_x\Delta w^1_{t_{i}} + h(x)\rhoi \Delta\tilde w_{t_{i}}+\alpha^5_i \right\}
|{\cal F}^{\tilde w}_{t_{i},T}\right],\hspace{3cm}
\end{eqnarray*}
where again, $\alpha_i^5 = o(\Delta t_i)$ in the same sense.

~

\noindent
{\bf 6}. Now, note that ${\cal F}^{\tilde w}_{t_{i-1},T}={\cal F}^{\tilde w}_{t_{i-1},t_{i}}
\bigvee
{\cal F}^{\tilde w}_{t_{i},T}$ and, moreover this $\sigma$-field is
independent from $w^1$.
Using the regular calculus for conditional expectations (cf. \cite{Ro}),
we get

\[
\tilde  \mathbb  E[V(t_{i},T,x)\rhoi |\Fw _{t_{i-1},T}]
=\tilde  \mathbb  E[V(t_{i},T,x) \rhoi |\Fw_{t_{i},T}],
\]
and in the same manner we can replace $\sigma$-fields $\Fw_{t_{i-1},T}$
by $\Fw_{t_{i},T}$ in all expressions in the previous step. 
Also, $\tilde  \mathbb  E\left[\Delta w^1_{t_{i}}|{\cal F}^{\tilde w}_{t_{i-1},T}\right] = 0$ due to the independence of $w^1$ and $\tilde w$ with respect to the measure $\tilde  \mathbb  P$, and $(\Delta w^1_{t_{i}})^2 \approx \Delta t_i$. 
Hence, we obtain

\new{\X}{X(t_i,T,x)}

\begin{eqnarray*}
&\tilde  \mathbb  E[V(t_{i},T,x)\rho^{-1}_{t_{i-1},T}|{\cal F}^{\tilde w}_{t_{i-1},T}]
\\\\
&=\tilde  \mathbb  E[V(t_{i},T,x)\rhoi |\Fw_{t_{i-1},T}]
\\\\
&\displaystyle + \tilde  \mathbb  E \left[\frac12  V(t_{i},T,x)(\rhoi)_{xx}
+ V_x(t_{i},T,x)(\rhoi)_x+\frac12 V_{xx}(t_{i},T,x)\rhoi |\Fw_{t_{i-1},T} \right]\Delta t_i
\\\\
&+\tilde  \mathbb  E[ f(x)(V_x(t_{i},T,x)\rhoi + V(t_{i},T,x)(\rhoi)_x)|\Fw_{t_{i-1},T}] \Delta t_i
\\\\
&+\tilde  \mathbb  E[V(t_{i},T,x) h(x)\rhoi \Delta\tilde w_{t_{i}}|\Fw_{t_{i-1},T}]
+ \alpha^6_i
\\\\
&\displaystyle =\tilde  \mathbb  E[V(t_{i},T,x)\rhoi |\Fw_{t_{i}}] 
\Delta t_i
+\tilde  \mathbb  E \left[\frac12 (V(t_{i},T,x)\rhoi)_{xx}|\Fw_{t_{i},T}\right] \Delta t_i
\\\\
&+\tilde  \mathbb  E[f(x)(V(t_{i},T,x)\rhoi)_x|\Fw_{t_{i},T}] \Delta t_i
+\Delta\tilde w_{t_{i}} \tilde  \mathbb  E[h(x)V(t_{i},T,x)\rhoi|\Fw_{t_{i},T}]
+\alpha^6_i
\\\\
&\displaystyle =v^g(t_i,x)+\frac12 v^g_{xx}(t_i,x)\Delta t_i+f(x)v^g_x(t_i,x)\Delta t_i
+ h(x)v^g(t_i,x)\Delta\tilde w_{t_{i}}+\alpha^6_i,
\end{eqnarray*}
with a similar property for $\alpha^6_i$: $\alpha^6_i = o(\Delta t_i)$ in the mean square sense. The last equality in this calculus  holds true
because of the possibility to change the order of integration
and derivation with respect to the $x$ variable.

~

\noindent
{\bf 7}. Therefore, we obtain the equality

\begin{eqnarray*}
&v^g(t_0,x)-v^g(T,x)
\\\\
& \displaystyle =\sum_i \left\{\frac12 v^g_{xx}(t_i,x)+f(x)v^g_x(t_i,x)\right\}\Delta t_i
+\sum_i h(x)v^g(t_i)\Delta\tilde w_{t_{i}} + \alpha^7,
\end{eqnarray*}
with $\alpha^7 =o(1)$ in the mean square sense as $\sup_i\Delta t_i\to 0$. Letting $\sup_i\Delta t_i\to 0$, we get from here the desired integral equality
(\ref{eq89}). The Theorem \ref{Th1} is proved.

\section*{Acknowledgement}
This work was prepared within the framework of a subsidy granted to the HSE by the Government of the Russian Federation for the implementation of the Global Competitiveness Program.


\begin{thebibliography}{99}

\bibitem{BlFr} 
Yu. N. Blagoveshchenskii, M. I. Freidlin,  
Certain properties of diffusion processes depending on a parameter. Math. USSR Doklady, 2(3) (1961), 633-636.





\bibitem{Kry-MP}
N. V. Krylov, On the selection of a Markov process from a system of processes and the construction of quasi-diffusion process, Math. USSR Izv., 7(3) (1973), 691-709. 


\bibitem{Kry} N. V. Krylov, 
On explicit formulas for solutions of evolutionary SPDE's 
(a kind of introduction to the theory). 
Lecture Notes in Control and Information Sciences, 1992, vol. 176, p. 153-164.


\bibitem{KR}
N. V. Krylov, B. L. Rozovsky, Stochastic differential equations and diffusion processes, Russian Math. Surveys, 37(6) (1982), 81-105.






\bibitem{KZ}
N. V. Krylov, A. Zatezalo, 
A Direct Approach to Deriving Filtering Equations for Diffusion Processes, 
Applied Mathematics and Optimization, 42(3) (2000), 315-332.




\bibitem{Kunita}
H. Kunita, Stochastic differential equations and stochastic flows of diffeomorphisms, 
\'Ecole d'\'Et\'e de Probabilit\'es de Saint-Flour XII (1982), 
Lecture Notes in Mathematics, vol. 1097 (1984), 143-303.





\bibitem{P}
E. Pardoux, 
Stochastic partial differential equations and filtering of diffusion processes, 
Stochastics, 3(1-4) (1980), 127-167.
 




\bibitem{Ro} 
B. L. Rozovski, Stochastic Evolution Systems. 
Kluwer Acad. Publ.: Dordreht et al., 1990. 






\bibitem{Ve} 
A. Yu. Veretennikov, Inverse Diffusion and Direct Derivation of Stochastic Liouville Equations. Soviet Math. Notes,  33(5-6) (1983), 397-400. 

\bibitem{Ve95}
A. Yu. Veretennikov,  
On backward filtering equations for SDE systems (direct approach). In: Stochastic Partial Differential Equations, ed. by A. Etheridge. London Math. Soc. Lecture Notes Series, Cambridge Univ. Press, vol. 216 (1995), 304-311.





\end{thebibliography}
\end{document}